\documentclass[preprint,12pt]{elsarticle}

\usepackage{amsmath, amssymb, amsfonts, amsthm, amscd}
\makeatletter
\renewcommand*\env@matrix[1][*\c@MaxMatrixCols c]{%
  \hskip -\arraycolsep
  \let\@ifnextchar\new@ifnextchar
  \array{#1}}
\makeatother
\usepackage{manfnt}
\usepackage{mathtools}
\usepackage{comment}
\usepackage{float}
\usepackage{fullpage}
\usepackage{url}
\usepackage[all]{xy}
   \SelectTips{cm}{10}
\usepackage{setspace}
\usepackage{txfonts}
\usepackage{enumitem}
\usepackage{verbatim} 
\numberwithin{equation}{section}
\newtheorem{theorem}{Theorem}[section]

\newtheorem{lemma}[theorem]{Lemma}
\newtheorem{prop}[theorem]{Proposition}
\newtheorem{defn}[theorem]{Definition}
\newtheorem{conjecture}[theorem]{Conjecture}

\newtheorem{cor}[theorem]{Corollary}

\theoremstyle{remark}
\newtheorem{rmk}{Remark}[theorem]
\newtheorem*{rmk*}{Remark}

\newcommand*\mathmeta[1]{\DOTSB\mathchoice{\;#1\;}{#1}{#1}{#1}}

\newcommand*\renewmetasymbol[2]{\renewcommand#1{\DOTSB\mathmeta{#2}}}

\renewmetasymbol\implies{\Rightarrow}
\renewmetasymbol\impliedby{\Leftarrow}
\renewmetasymbol\iff{\Leftrightarrow}

\newcommand{\resultall}{43/246}

\newcommand{\Z}{\mathbb{Z}}

\newcommand{\Q}{\mathbb{Q}}

\newcommand{\R}{\mathbb{R}}

\newcommand{\e}{\epsilon}

\DeclareMathOperator{\ord}{ord}
\numberwithin{equation}{section}
\journal{Acta Arithmetica}

\begin{document}

\begin{frontmatter}
\title{Heights of Rational Points on Mordell Curves}
\author[inst1]{Alan Zhao}
\affiliation[inst1]{organization={Columbia University Department of Mathematics},
            addressline={2990 Broadway, Room 509, MC 4406}, 
            city={New York},
            state={NY},
            postcode={10027}, 
            country={USA}}

\begin{abstract}
We conjecture a lower bound for the minimal canonical height of non-torsion rational points on a natural density $1$ subset of the sextic twist family of Mordell curves. We then establish a lower bound that yields a partial result towards this conjecture.
\end{abstract}


\end{frontmatter}

\tableofcontents
\section{Introduction}
\subsection{Summary}
In \cite{le2016height}\cite{boudec2018height}, Le Boudec establishes a lower bound for the minimal canonical height of non-torsion rational points on a natural density $1$ subset of a generic quadratic twist family $E^{(\delta )}(A,B): \delta y^2=x^3+Ax+B$, where $A, B \in \Z$ such that $4A^3 + 27B^2 \neq 0$ and $\delta \in \Z_{>0}$ square-free. Stronger bounds are proven upon specializing to the case of $A = -1$ and $B = 0$. To formally state these results, we will need to set some notation. Let $\mathcal{B} \in \R_{> 0}$ and set $S_2(\mathcal{B}) = \{\delta  \in \Z_{>0}: \delta \leq \mathcal{B}, \delta \text{ square-free}\}$. Given an elliptic curve $E/\Q$, let $\hat{h}_E$ be the canonical height. Now, we can make

\begin{defn} \label{correspond}
Let $\log \eta_{\delta}(A,B)= \min \{ \hat{h}_{E^{(\delta )}(A,B)}(P): P \in E^{(\delta )}(A,B)(\Q) \setminus E^{(\delta )}(A,B)(\Q)_{\mathrm{tors}} \}$. If $E^{(\delta )}(A,B)(\Q) \setminus E^{(\delta )}(A,B)(\Q)_{\mathrm{tors}}$ is empty, then $\eta_{\delta}(A,B) = \infty$.
\end{defn}

The works \cite{le2016height}\cite{boudec2018height} establish lower bounds for $\eta_{\delta}(A,B)$ and $\eta_{\delta}(-1, 0)$. The motivation for doing so comes from recent literature studying the analogy between number fields and elliptic curves, a discussion of which may be found in \cite[\S 1.2]{le2016height}. The upshot of this discussion is

\begin{conjecture}[{\cite[Conjecture A]{le2016height}}] \label{pierre}
Let $\e>0$ and $A,B \in \Z$ such that $4A^3+27B^2 \neq 0$. The set of $\delta  \in S_2(\mathcal{B})$ such that
\begin{equation*}
    \eta_{\delta }(A,B) > e^{\delta ^{1/2-\e}}
\end{equation*}
has natural density $1$ in $S_2(\mathcal{B})$ as $\mathcal{B} \to \infty$.
\end{conjecture}

The list of results in this direction is given in
\begin{theorem}[{\cite[Theorem 1]{le2016height}\cite[Theorem 2]{le2016height}\cite[Theorem 1]{boudec2018height}}] \label{hard}
Let $\e > 0$ and $A,B \in \Z$ such that $4A^3+27B^2 \neq 0$. The sets of $\delta  \in S_2(\mathcal{B})$ such that 
\begin{itemize}
    \item $\eta_{\delta }(A,B) > \delta ^{1/4 - \e}$ 
    \item $\eta_{\delta }(-1,0) > \delta ^{5/8 - \e}$ 
    \item $\eta_{\delta }(-1,0) > \delta ^{0.845}$ 
\end{itemize}
have natural density $1$ in $S_2(\mathcal{B})$ as $\mathcal{B} \to \infty$ in the first two cases and positive natural density in $S_2(\mathcal{B})$ as $\mathcal{B} \to \infty$ in the third case.
\end{theorem}

We now seek to expand these results to the sextic twist family of Mordell curves
\[E_d: y^2z=x^3+dz^3,\]
where $d \in \Z$ sixth-power free.

\subsection{Bounds for Heights of Points on $E_d$}
In analogy to Definition \ref{correspond}, we make

\begin{defn}
Let $\log \zeta_d= \min \{ \hat{h}_{E_d}(P): P \in E_d(\Q) \setminus E_d(\Q)_{\mathrm{tors}} \}$. If $E_d(\Q) \setminus E_d(\Q)_{\mathrm{tors}}$ is empty, then $\zeta_d= \infty$.
\end{defn}

Let $S_6(\mathcal{B}) = \{d \in \Z: |d| \leq \mathcal{B}, \text{ } d \text{ sixth-power free}\}$. Before we begin, let us set up a target analogous to Conjecture \ref{pierre}. 

\begin{conjecture} \label{myconjecture}
Let $\e > 0$. The set of $d \in S_6(\mathcal{B})$ such that 
\begin{equation*}
\zeta_d > e^{|d|^{1/6-\e}}
\end{equation*}
has natural density $1$ in $S_6(\mathcal{B})$ as $\mathcal{B} \to \infty$.
\end{conjecture}
We choose the exponent of $|d|$ in Conjecture \ref{myconjecture} so that it specializes to Conjecture \ref{pierre}. More precisely, when $d > 0$ (resp. $d < 0$), we have that $E_{d^3}$ is isomorphic to $E^{(d)}(0, 1)$ (resp. $E^{(|d|)}(0, -1)$ ) by the $\Q$-linear change of variables $x \mapsto |d|x$ and $y \mapsto d^2y$. So, by replacing $d$ with $d^3$ in Conjecture \ref{myconjecture}, we obtain the statement of Conjecture \ref{pierre} for $A = 0$ and $B = \pm 1$.

Now, it is important to note that $E^{(\delta)}(A,B)$ and $E_d$ are conjectured to have the same rank distributions. In particular, Goldfeld's conjecture \cite[Conjecture B]{goldfeld1979conjectures} predicts that the set of $\delta \in S_2(\mathcal{B})$ such that $E^{(\delta)}(A,B)$ has analytic rank $0$ (resp. $1$) has natural density $1/2$ (resp. $1/2$) in $S_2(\mathcal{B})$ as $\mathcal{B} \to \infty$. We then have the same distribution for algebraic ranks $0$ and $1$ by the work of Gross--Zagier \cite{gross1986heegner} and Kolyvagin \cite{kolyvagin1989finiteness} on the Birch and Swinnerton-Dyer conjecture. Because of our definition that $\eta_{\delta}(A,B) = \infty$ for $E^{(\delta)}(A,B)$ of rank $0$, assuming Goldfeld's conjecture, it would suffice to consider $E^{(\delta)}(A,B)$ of rank $1$ to prove Conjecture \ref{pierre}. 

A similar reduction is expected in the case of the sextic twist family of Mordell curves. It is widely believed that Goldfeld's conjecture, with $E^{(\delta)}(A,B)$ replaced by $E_d$, should also hold. Formally, it is expected that the set of $d \in S_6(\mathcal{B})$ such that $E_d$ has analytic rank $0$ (resp. $1$) has natural density $1/2$ (resp. $1/2$) as $\mathcal{B} \to \infty$. Thus, it would also suffice to consider $E_d$ of rank $1$ to prove Conjecture \ref{myconjecture}.

Furthermore, the progress made towards Goldfeld's conjecture in \cite[Theorem 1.5]{smith2016congruent} is the key ingredient of the proof for the third and strongest bound in Theorem \ref{hard}. Similar progress has been made for Goldfeld's conjecture for the sextic twist family of Mordell curves, the details of which may be found in \cite[Theorem 1.8]{kriz2019goldfeld}. A lower bound for $\zeta_d$ on a natural density $1$ subset of this family is given in Theorem \ref{20} below, which is the main result of this paper. So, by possibly restricting ourselves to a positive natural density subset, we can also expect to be able to strengthen the bound of this theorem. This would yield improvement analogous to that which comes from strengthening the second bound to the third bound in Theorem \ref{hard}.

Not making this restriction for now, we set out to make progress towards Conjecture \ref{myconjecture} by proving
\begin{theorem}[Corollary \ref{200}] \label{20}
Let $\e > 0$. The set of $d \in S_6(\mathcal{B})$ such that 
\begin{equation*}
    \zeta_d > |d|^{\resultall-\e}
\end{equation*}
has natural density $1$ in $S_6(\mathcal{B})$ as $\mathcal{B} \to \infty$.
\end{theorem}

We will now proceed as follows: in \S 2, we provide a parameterization of the rational points on $E_d$ and then prove the sharp lower bound $\zeta_d \gg |d|^{1/36}$. In \S 3, we tighten this bound under the assumption that the square part of $d$ is small. This will be sufficient to prove Theorem \ref{20}.

\begin{rmk*}
The curve $E_d$ has a rational $3$-isogeny. Thus, it would be interesting to see if an application of $3$-descent \cite[Theorem 3.1]{cohen2009elementary}\cite[Theorem 4.1]{cohen2009elementary} may also be used to strengthen Theorem \ref{20}.
\end{rmk*}

\section*{Acknowledgements} 
I would like to thank Chao Li of Columbia University for his invaluable mentorship, and especially for his willingness to help me remotely over numerous video calls during the pandemic. This research has been generously supported through a grant funded by the I.I. Rabi Scholars Program at Columbia University. 

\section{Preliminaries}
Because of the expected sharp bound of $\zeta_d \gg |d|^{1/36}$, we let $\alpha \in \R_{> 0}$ and work with the quantity
\begin{align} \label{lookbackto}
    &N_{\alpha}(\mathcal{B}) := \#\{d \in S_6(\mathcal{B}): \zeta_d \leq |d|^{1/36 + \alpha}\}.
\end{align}
This quantity tells us how many $d \in S_6(\mathcal{B})$ fail to satisfy $\zeta_d > |d|^{1/36+\alpha}$. So, if $N_{\alpha}(\mathcal{B}) = o(\mathcal{B})$, then, since $\#S_6(\mathcal{B}) \gg \mathcal{B}$, the set of $d \in S_6(\mathcal{B})$ such that $\zeta_d > |d|^{1/36+\alpha}$ has natural density $1$ in $S_6(\mathcal{B})$ as $\mathcal{B} \to \infty$. The goal will now be to maximize $\alpha$ with respect to the condition that $N_{\alpha}(\mathcal{B}) = o(\mathcal{B})$.

\subsection{Reducing from Canonical Height to Logarithmic Height}
Let $h: \mathbb{P}^1(\overline{\Q}) \to \R_{\geq 0}$ be the logarithmic absolute height. For points $[a:b]$ defined over $\Q$, we have that $h([a:b]) = \log \max(|a|, |b|)$, where the representative $[a:b]$ is chosen so that $a,b \in \Z$ and $\gcd(a,b) = 1$. Given an elliptic curve $E/\Q$, we may consider any $f \in \Q(E)$ as a function $E \to \mathbb{P}^1$ given by $P \mapsto [f(P): 1]$, which we also denote as $f$. Now, let $h_f(P) = h(f(P))$.
The quantity $h_f$ is related to the quantity $\hat{h}_E$ by
\begin{lemma}[{\cite[\S VIII.9, Theorem 9.3(e)]{silverman2009arithmetic}}] \label{pseudo}
Let $E/\Q$ be an elliptic curve, $P \in E(\overline{\Q})$, and $f \in \Q(E)$ even. Then,
\begin{equation*}
    \hat{h}_E(P) = \frac{1}{\deg f}h_f(P) + O_{E,f}(1).
\end{equation*}
\end{lemma}

We will use this relation to give an approximation of the canonical height in
\begin{lemma}
Let $P_1 \in E_1(\overline{\Q})$ and $f \in \Q(E_1)$ be given by $f(P_1) = x(P_1)^3/y(P_1)^2$, where $x$ and $y$ are Weierstrass coordinates for the curve $E_1$. For any $d \in \Z$ sixth-power free and $P \in E_d(\overline{\Q})$, we have
\begin{equation*}
    \hat{h}_{E_d}(P) = \frac{1}{6}h_f(P) + O(1).
\end{equation*}
\end{lemma}
\begin{proof}
Begin by noting that $f$ is even since $x$ and $y^2$ are both even. Also note that we have an isomorphism $\kappa: E_d \to E_1$ over $\overline{\Q}$ given by $\kappa([x:y:z]) = [q_3x: q_2y:z]$, where $q_2, q_3 \in \overline{\Q}$ such that $q_2^2 = q_3^3 = d^{-1}$. Because $\hat{h}_{E_1}$ is invariant under $\overline{\Q}$-isomorphisms of $E_1$, by applying Lemma \ref{pseudo} we have that
\begin{equation*}
    \hat{h}_{E_d}(P) = \hat{h}_{E_1}(\kappa(P)) = \frac{1}{\deg f}h_f(\kappa(P)) + O(1).
\end{equation*}
Note that the $O(1)$ term no longer has dependencies as it does in Lemma \ref{pseudo} because we have made a choice of even function $f$ and an elliptic curve $E_1$. Since $h_f(\kappa(P)) = h_f(P)$, it now suffices to prove that $\deg f= 6$. 

To begin, observe that $x^3$ and $y^2$ have the same order pole at the base point of $E_1$, no other poles, and do not share any zeroes. Thus, the poles of $f$ are precisely the zeroes of $y$, and thus
\begin{equation} \label{sum}
    \deg f = \sum_{\ord_Q(f) < 0} -\ord_Q(f) = \sum_{\ord_Q(y) > 0} \ord_Q(y^2).
\end{equation}
Then, by the theory of divisors,
\begin{equation} \label{sum2}
    \sum_{\ord_Q(y) > 0} \ord_Q(y^2) = \sum_{\ord_Q(y) < 0} \ord_Q(y^2) = 6,
\end{equation}
where the second equality is true since $y$ has only one pole of order $3$. Combining \eqref{sum} and \eqref{sum2}, we have that $\deg f = 6$, as wished.
\end{proof}
For $P \in E_d(\overline{\Q})$, the lemma implies that
\begin{equation*}
e^{h_{x^3/y^2}(P)} \ll e^{6\hat{h}_{E_d}(P)} \ll e^{h_{x^3/y^2}(P)}.
\end{equation*}
This yields
\begin{cor} \label{6}
For any $a \geq 0$ and $P \in E_d(\overline{\Q})$, we have that $e^{{h}_{x^3/y^2}(P)} \gg |d|^{a}$ (resp. $e^{{h}_{x^3/y^2}(P)} \ll|d|^{a}$) if and only if $e^{\hat{h}_{E_d}(P)} \gg |d|^{a/6}$ (resp. $e^{\hat{h}_{E_d}(P)} \ll |d|^{a/6}$).
\end{cor}
\subsection{Parameterization Data of $E_d$} 
To use the function $h_{x^3/y^2}$, it will be necessary to extract the common factors of $x$ and $y$, which is the goal of this section. To begin, define $\square: \Z_{\neq 0} \to \Z_{>0}$ such that $\square(n)$ is the positive square-part of $n$. Then, we may state

\begin{prop} \label{?impor}
Let $\e > 0$. Set
\begin{equation} \label{G}
   G_{\mathcal{B}, \alpha, \e} :=  \left\{ (D,X,Y,Z) \in \Z^3 \times \Z_{>0}: \begin{aligned} & D \in S_6(\mathcal{B}) \\ & \gcd(X,Y,Z) = 1 \\ & Y^2Z=X^3+DZ^3 \\ &e^{h_{x^3/y^2}([X:Y:Z])} \ll \mathcal{B}^{1/6+6\alpha} \\ &\square(D) \leq \mathcal{B}^{\e} \end{aligned}\right\}
\end{equation}
and
\begin{equation} \label{paramm}
    H_{\mathcal{B}, \alpha, \e} := \left \{ (C, B_1, B_2, D_2, Y_2, X_1) \in \Z_{>0}^3 \times \Z_{\neq 0} \times \Z^2: \begin{aligned} & C^3B_2^2D_2 \in S_6(\mathcal{B})\\ & \gcd(B_1,CX_1B_2Y_2) = 1 \\ & \gcd(B_2X_1, CY_2) = 1 \\ & CY_2^2 = B_2X_1^3 + D_2B_1^6 \\ &|B_2X_1^3|, |CY_2^2| \ll \mathcal{B}^{1/6+6\alpha} \\ & |C|, |B_2| \leq \mathcal{B}^{\e}\end{aligned}\right \}.
\end{equation}
Then, there are injections of sets $g_{\mathcal{B}, \alpha, \e}: G_{\mathcal{B}, \alpha, \e} \hookrightarrow  H_{\mathcal{B}, \alpha, \e}$ such that, given $(d, x, y, z) \in G_{\mathcal{B}, \alpha, \e}$, the data of its image $(c,b_1,b_2,d_2,y_2,x_1) \in H_{\mathcal{B}, \alpha, \e}$ gives a parameterization $d = c^3b_2^2d_2$, $x = cb_1b_2x_1$, $y=c^2b_2y_2$, and $z = b_1^3$. When $\e \geq 1/2$, $g_{\mathcal{B}, \alpha, \e}$ is an isomorphism.
\end{prop}
\begin{proof}
We begin by constructing a map $g_{\mathcal{B}, \alpha, \e}: G_{\mathcal{B}, \alpha, \e} \to H_{\mathcal{B}, \alpha, \e}$ that gives the required parameterization. In what follows, given $(d,x,y,z) \in G_{\mathcal{B}, \alpha, \e}$, we will define a tuple $(c,b_1,b_2,d_2,y_2,x_1) \in H_{\mathcal{B}, \alpha, \e}$ and set $g_{\mathcal{B}, \alpha, \e}(d,x,y,z) = (c,b_1,b_2,d_2,y_2,x_1)$. Then, we will prove the remaining desired properties of $g_{\mathcal{B}, \alpha, \e}$.

Since $z > 0$, at least one of $x$ and $y$ is non-zero. Thus, we may let $b_0 = \gcd(x,y)$ and write $x = b_0x_0$ and $y = b_0y_1$ such that $\gcd(x_0,y_1) = 1$. Substitute these expressions into the relation $y^2z = x^3+dz^3$ to obtain
\begin{equation} \label{eq1}
    b_0^2y_1^2z = b_0^3x_0^3 + dz^3.
\end{equation}
This tells us that $b_0^2 \mid dz^3$. Since $\gcd(x,y,z) = 1$, we have that $\gcd(b_0,z) = 1$, and thus $b_0^2 \mid d$. Write $d = b_0^2d_1$, and substitute this into \eqref{eq1} to obtain
\begin{equation} \label{A}
    y_1^2z =b_0x_0^3 + d_1z^3.
\end{equation}

Now, let $b_1 = \gcd(x_0, z)$. Then, $b_1^3 \mid y_1^2z$. Since $\gcd(x,y,z) = 1$, we have that $\gcd(b_1,y_1) = 1$, and thus $b_1^3 \mid z$. We may then write $x_0 = b_1x_1$ and $z = ub_1^3$ such that $\gcd(x_1, ub_1) = 1$. Since $z, b_1 > 0$, to prove that $u = 1$ it will suffice for us to prove that $u$ divides a power of $x_1$. We begin by substituting the new expressions for $x_0$ and $z$ into \eqref{A} to obtain
\begin{equation} \label{eq2}
    y_1^2u = b_0x_1^3+d_1u^3b_1^6.
\end{equation}
This tells us that $u \mid b_0x_1^3$. Since $\gcd(b_0, z) = 1$, we have that $\gcd(b_0, u) = 1$. Thus, $u \mid x_1^3$ as wished, and so $u = 1$ and \eqref{eq2} simplifies to
\begin{equation} \label{eq3}
    y_1^2 = b_0x_1^3+d_1b_1^6.
\end{equation} 

Finally, let $c = \gcd(b_0, y_1)$ and write $b_0 = cb_2$ and $y_1 = cy_2$ such that $\gcd(b_2,y_2) = 1$. Substituting these expressions into \eqref{eq3}, we obtain 
\begin{equation} \label{eq4}
    c^2y_2^2 = cb_2x_1^3 + d_1b_1^6,
\end{equation} 
which implies that $c \mid d_1b_1^6$. Since $\gcd(b_0, z) = 1$, we have that $\gcd(b_0, b_1) = 1$. By construction, $c \mid b_0$, and so we have that $\gcd(c,b_1) = 1$. Thus, $c \mid d_1$, and we may now write $d_1 = cd_2$. Substituting this into \eqref{eq4}, we obtain 
\begin{equation} \label{eq5}
cy_2^2 = b_2x_1^3 + d_2b_1^6.
\end{equation}
This tells us that $\gcd(c,b_2) \mid d_2b_1^6$, and since $\gcd(c,b_1) = 1$, we have that $\gcd(c,b_2) \mid d_2$. But $c^3b_2^2d_2$ is sixth-power free, and so we must have that $\gcd(c,b_2) = 1$.

Tracing through the variable declarations thus far, we can see that we have defined a tuple $(c,b_1,b_2,d_2,y_2,x_1)$ that gives parameterizations of $d,x,y$, and $z$ as in the lemma statement. It now remains to show that $(c,b_1,b_2,d_2,y_2,x_1) \in H_{\mathcal{B}, \alpha, \e}$, which we will do by directly checking that it satisfies all conditions in \eqref{paramm}.

The GCD is always positive, and so $c, b_1, b_2 > 0$. From the parameterization of $d$, we obtain the first condition of \eqref{paramm}, and thus we also have that $d_2 \neq 0$. From the content directly after \eqref{eq5}, we have that $\gcd(c,b_2) = 1$. From the content just before
\begin{itemize}
    \item \eqref{eq1}, we have that $\gcd(x_0, y_1) = \gcd(b_1x_1,cy_2) = 1$.
    \item \eqref{A}, we have that $\gcd(b_0,z) = \gcd(cb_2,b_1^3) = 1$, and hence $\gcd(b_2, b_1) = 1$.
    \item \eqref{eq2}, we have that $\gcd(x_1,ub_1) = \gcd(x_1,b_1) = 1$ since $u = 1$.
    \item \eqref{eq4}, we have that $\gcd(b_2,y_2) = 1$.
\end{itemize}
Thus, the second and third conditions of \eqref{paramm} are satisfied. The fourth condition of \eqref{paramm} holds by \eqref{eq5}. Given the second and third condtions of \eqref{paramm}, the fourth condition of \eqref{G}, and the parameterizations of $d,x,y$, and $z$, we have that 
\begin{equation}\label{forlater}
e^{h_{x^3/y^2}([x:y:z])} = e^{h\left( \left[\left(\frac{cb_2b_1x_1}{b_1^3} \right)^3: \left( \frac{c^2b_2y_2}{b_1^3} \right)^2\right]\right)} =  \max\{|b_2x_1^3|, |cy_2^2|\} \ll \mathcal{B}^{1/6+6\alpha}.
\end{equation}
This gives the fifth condition of \eqref{paramm}. Finally, since $d = c^3b_2^2d_2$, we have that $\square(d) \geq cb_2$. Since $\square(d) \leq \mathcal{B}^{\e}$, we obtain $|cb_2| \leq \mathcal{B}^{\e}$. This gives the sixth condition of \eqref{paramm}. In sum, $(c,b_1,b_2,d_2,y_2,x_1)$ satisfies all conditions in \eqref{paramm}. This completes the construction of $g_{\mathcal{B}, \alpha, \e}$.

The fact that $g_{\mathcal{B}, \alpha, \e}$ is injective is immediate: for any two tuples $(d_i, x_i, y_i, z_i) \in G_{\mathcal{B}, \alpha, \e}$ (where $i \in \{1,2\}$) which map to the same tuple $(c, b_1, b_2, d_2,y_2, x_1) \in H_{\mathcal{B}, \alpha, \e}$, from the parameterization afforded by $g_{\mathcal{B}, \alpha, \e}$, we obtain that $(d_1,x_1,y_1,z_1) = (d_2,x_2,y_2,z_2)$.

For the rest of the proof, let $\e \geq 1/2$. It remains to show that $g_{\mathcal{B}, \alpha, \e}$ is a surjection. We show that it has a right inverse $\hat{g}_{\mathcal{B}, \alpha, \e}: H_{\mathcal{B}, \alpha, \e} \to G_{\mathcal{B}, \alpha, \e}$ given by $(c,b_1,b_2,d_2,y_2,x_1) \mapsto (d,x,y,z) = (c^3b_2^2d_2, cb_1b_2x_1, c^2b_2y_2, b_1^3)$. To see that this map is well-defined, we show that $(d,x,y,z)$ satisfies each condition in \eqref{G}. Taking the conditions from \eqref{paramm} and using
\begin{itemize}
    \item that $b_1 > 0$, we have that $z > 0$.
    \item its first, we know that $d = c^3b_2^2d_2 \in S_6(\mathcal{B})$. 
    \item its second and third, we know that $\gcd(x,y,z) = \gcd(cb_1b_2x_1, c^2b_2y_2, b_1^3) = \gcd(cb_2, b_1^3) = 1$. 
    \item its fourth, we have that $y^2z = x^3+dz^3$.
    \item its second, third, and fifth, we obtain \eqref{forlater}.
    \item its first, $|d| \leq \mathcal{B}$, and thus $\square(d) \leq \mathcal{B}^{1/2}$.
\end{itemize}
Thus, all conditions of \eqref{G} are satisfied, and so $(d,x,y,z) \in G_{\mathcal{B}, \alpha, \e}$.  To show that $\hat{g}_{\mathcal{B}, \alpha, \e}$ is a right inverse, it suffices to show that if $g_{\mathcal{B}, \alpha, \e}(d,x,y,z) = (c', b_1', b_2', d_2', y_2', x_1')$, then we have that $(c',b_1', b_2',d_2', y_2', x_1') = (c,b_1,b_2,d_2,y_2,x_1)$.

From the parameterization that $g_{\mathcal{B}, \alpha, \e}$ provides, we have that $d = c'^3b_2'^2d_2'$, $x = c'b_1'b_2'x_1'$, $y = c'^2b_2'y_2'$, and $z = b_1'^3$. From the parameterization for $z$, we have that $b_1 = b_1'$. By the second and third conditions in \eqref{paramm}, we have that $\gcd(x,y) = cb_2 = c'b_2'$,
\begin{itemize}
    \item which implies that $cy_2 = \frac{y}{\gcd(x,y)} = c'y_2'$. Combining these equalities with the third condition in \eqref{paramm}, we have that $c = \gcd(cb_2, cy_2) = \gcd(c'b_2',c'y_2') = c'$, and hence $b_2 = b_2'$ and $y_2 = y_2'$.
    \item which implies that $b_1x_1 = \frac{x}{\gcd(x,y)} = b_1'x_1'$. Since $b_1 = b_1'$, we have that $x_1 = x_1'$. 
\end{itemize}
 By the parameterization for $d$, since $c = c'$ and $b_2 = b_2'$, it then follows that $d_2 = d_2'$. Gathering all equalities, we have shown that $(c,b_1,b_2,d_2,y_2,x_1) = (c', b_1', b_2', d_2',y_2',x_1')$, as wished.
\end{proof}
\begin{rmk} \label{tech}
Non-torsion rational points $P$ lying on curves in the sextic twist family of Mordell curves are in 1-to-1 correspondence with tuples $(d,x,y,z)$ such that $P \in E_d(\Q)$ is represented by $[x:y:z] \in \mathbb{P}^2(\Z)$. This representation is unique (and hence so is the tuple $(d,x,y,z)$ characterizing $P$) when we require that $\gcd(x,y,z) = 1$ and $z > 0$.
\end{rmk}
With the common divisors of $x$ and $y$ now examined, we can apply the results obtained in \S 2.1 to prove
\begin{prop}
For any $d \in \Z$ sixth-power free, we have that
\begin{equation*}
    \zeta_d \gg |d|^{1/36}.
\end{equation*}
Moreover, this bound is sharp.
\end{prop}
\begin{proof}
For $P \in E_d(\Q) \setminus E_d(\Q)_{\mathrm{tors}}$, take its unique representative $[x:y:z] \in \mathbb{P}^2(\Z)$ satisfying the conditions discussed in Remark \ref{tech}. In the notation of Proposition \ref{?impor}, we then have that $(d,x,y,z) \in G_{\mathcal{B}, \alpha, 1/2}$ for $\mathcal{B}$ sufficiently large. We may then apply the isomorphism $g_{\mathcal{B}, \alpha, 1/2}$ in Proposition \ref{?impor} to $(d,x,y,z)$ to parameterize $d,x,y$, and $z$. Then, by the equalities in \eqref{forlater} and comparing terms in \eqref{eq5}, we obtain
\begin{equation*}
    e^{h_{x^3/y^2}(P)} = \max \{|b_2x_1^3|, |cy_2^2|\} \gg |d_2b_1^6|.
\end{equation*}
We can now observe that
\begin{align} 
    \max \{|b_2x_1^3|, |cy_2^2|\} & \label{align}\gg |cy_2^2|^{1/2}|b_2x_1^3|^{1/3}|d_2b_1^6|^{1/6} \\ &\label{align2}\gg |c|^{1/2}|b_2|^{1/3}|d_2|^{1/6} = |d|^{1/6}.
\end{align}
By Corollary \ref{6}, this implies that $e^{\hat{h}_{E_d}(P)}\gg |d|^{1/36}$, and hence $\zeta_d \gg |d|^{1/36}$. It now remains to prove that this inequality for $\zeta_d$ is sharp. To do this, we must show that the lower bounds for \eqref{align} and \eqref{align2} are attained for infinitely many $d$.

Taking $b_2 \in \Z_{>0}$, let $ (c,b_1,b_2,d_2,y_2,x_1) = T_{b_2} := (2b_2+1, 1, b_2, b_2+1, 1, 1)$. Then, the lower bound of \eqref{align2} is attained since $y_2 = x_1 = b_1 = 1$, and that of \eqref{align} is attained since for $n \in \{|b_2|, |c|, |d_2|\}$, we have that $n \ll \max\{ |b_2|, |c| \} \ll n$. Here, the implied constants depend on our choice that $(c,b_1,b_2,d_2,y_2,x_1) = T_{b_2}$.

Given $b_2$, take $\e = 1/2$ and choose $\mathcal{B}$ sufficiently large so that $T_{b_2}$ satisfies all conditions in \eqref{paramm} except for possibly the sixth-free part of the first condition. Thus, we have that $h(b_2)$ sixth-power fre implies $T_{b_2} \in H_{\mathcal{B}, \alpha, 1/2}$, where $h(b_2) = (2b_2+1)^3b_2^2(b_2 + 1)$. Applying the isomorphism $g_{\mathcal{B}, \alpha, 1/2}$ in Proposition \ref{?impor}, we obtain the tuple $(h(b_2), b_2(2b_2 + 1), b_2(2b_2+1)^2, 1) \in G_{\mathcal{B}, \alpha, 1/2}$ --- and thus, by Remark \ref{tech}, a rational point on $E_{h(b_2)}$ --- whose parameterization attains the lower bounds in \eqref{align} and \eqref{align2}. In sum, values $d = h(b_2)$ such that $\square(h(b_2)) = 1$ attain the lower bound of the inequality $\zeta_d \gg |d|^{1/36}$.

To prove that this inequality is sharp, it now suffices to show that there are infinitely many values $h(b_2)$ such that $h(b_2)$ is sixth-power free. Since $h(b_2)$ is strictly increasing when $b_2 > 0$, we may equivalently show that there are infinitely many $b_2$ such that $h(b_2)$ is sixth-power free. To do this, it will suffice to show that there are infinitely many $b_2$ such that $b_2, b_2+1$, and $2b_2+1$ are square-free. Since $b_2$, $b_2+1$, and $2b_2 + 1$ are pairwise coprime, it suffices to show that there are infinitely many $b_2$ such that $b_2(b_2+1)(2b_2+1)$ is square-free. To this end, let $A \in \R_{>0}$ and define
\begin{equation*}
    \mathcal{N}(A) :=  \#\left\{b_2: \begin{aligned} & b_2 \in \Z_{>0} \cap [1, A] \\ & b_2(b_2+1)(2b_2+1)\text{ square-free}\end{aligned}\right\}.
\end{equation*}
It remains to show that $\mathcal{N}(A) \to \infty$ as $A \to \infty$. But this follows from a direct application of \cite[Theorem 1.1]{booker2015square} to $\mathcal{N}(A)$, which tells us, for any $\e' > 0$, that $\mathcal{N}(A) \gg A^{1-\e'}$. Thus, by choosing $\e'$ sufficiently small, sharpness follows.
\end{proof} 

Having proven the promised sharp lower bound, the goal will now be to augment it with heavier technology. In particular, we will draw from estimates established using Heath-Brown's determinant method, and hence motivate the use of similar results in making progress towards Conjecture \ref{myconjecture}. 

Before we begin this process, we need to make some additional definitions. Let $\e > 0$, and define
\begin{equation*}
    N_{\alpha,\e}(\mathcal{B}) := \#\{d \in S_6(\mathcal{B}): \square(d) \leq \mathcal{B}^{\e}, \zeta_d \leq |d|^{1/36 + \alpha}\}
\end{equation*}
and
\begin{equation} \label{disc}
    N^*_{\alpha,\e}(\mathcal{B}) := \underset{\square(d) \leq \mathcal{B}^{\e}}{\sum_{d \in S_6(\mathcal{B})}}\#\{P \in E_d(\Q) \setminus E_d(\Q)_{\mathrm{tors}}: e^{\hat{h}_{E_d}(P)} \leq |d|^{1/36 + \alpha}\}.
\end{equation}
Because each $d$ counted by $N_{\alpha,\e}(\mathcal{B})$ corresponds to at least one $P \in E_d(\Q) \setminus E_d(\Q)_{\mathrm{tors}}$ with $e^{\hat{h}_{E_d}(P)} \leq |d|^{1/36 + \alpha}$, we have that $ N_{\alpha,\e}(\mathcal{B}) \leq  N^*_{\alpha,\e}(\mathcal{B})$.

With the notation in place, we can now turn towards the proof of Theorem \ref{20}.

\section{A Direct Count via Uniform Bounds}
In recent decades, great progress has been made on establishing bounds for the number of bounded integer zeroes of forms in $2$ or $3$ variables. We will use the following theorem of Heath-Brown to illustrate the utility of results like these in our analysis of the canonical height.
\begin{theorem}[{\cite[Theorem 3]{heath2002density}}] \label{bound}
Let $F \in \Z[x_1, x_2, x_3]$ be a singular form of degree $k$ irreducible over $\Q$. Let $N(F;\mathcal{B}_1, \mathcal{B}_2, \mathcal{B}_3)$, where $\mathcal{B}_i \geq 1$, be the number of integer zeroes of $F$ subject to the condition that $|x_i| \leq \mathcal{B}_i$. Finally, let $T = \max \mathcal{B}_1^{f_1}\mathcal{B}_2^{f_2}\mathcal{B}_3^{f_3}$, where the maximum is taken over all tuples $(f_1, f_2, f_3)$ such that the corresponding monomial $x_1^{f_1}x_2^{f_2}x_3^{f_3}$ occurs with non-zero coefficient.
Then, one has that
\begin{equation*}
    N(F;\mathcal{B}_1, \mathcal{B}_2, \mathcal{B}_3) \ll_{\e} T^{-k^{-2}}(\mathcal{B}_1\mathcal{B}_2\mathcal{B}_3)^{k^{-1}+\e}.
\end{equation*}
\end{theorem}

We apply this theorem to prove
\begin{prop} \label{step2}
Let $m(\alpha) = \min\left(1, \frac{1}{6} + 6\alpha\right)$ and $e(\alpha) = \e\left(m(\alpha) + \alpha + \frac{1}{36}\right)$. Then, one has that 
\begin{equation*}
N^*_{\alpha,\e}(\mathcal{B}) \ll_{\e} \mathcal{B}^{\frac{6}{49}(m(\alpha)+41\alpha)+ c_1+c_2 +2\e + e(\alpha)},
\end{equation*}
where $c_1 = \frac{1}{7 \cdot 6^2} - \frac{1}{6\cdot 7^2}  = \frac{1}{42^2} = \frac{1}{1764}$, and $c_2 = \frac{1}{12} + \frac{1}{18} = \frac{5}{36}$.
\end{prop}
\begin{proof}
By Corollary \ref{6}, we obtain 
\begin{equation*}
    N^*_{\alpha,\e}(\mathcal{B}) \ll \underset{\square(d) \leq \mathcal{B}^{\e}}{\sum_{d \in S_6(\mathcal{B})}}\#\{P \in E_d(\Q) \setminus E_d(\Q)_{\mathrm{tors}}:  e^{h_{x^3/y^2}(P)} \ll |d|^{1/6 + 6\alpha}\},
\end{equation*}
which implies that
\begin{equation} \label{p}
    N^*_{\alpha,\e}(\mathcal{B}) \ll \underset{\square(d) \leq \mathcal{B}^{\e}}{\sum_{d \in S_6(\mathcal{B})}}\#\{P \in E_d(\Q) \setminus E_d(\Q)_{\mathrm{tors}}: e^{h_{x^3/y^2}(P)} \ll \mathcal{B}^{1/6 + 6\alpha}\}.
\end{equation}
By Remark \ref{tech}, the set of points counted by the right-hand side of the above inequality may be identified with $G_{\mathcal{B}, \alpha, \e}$. Then, by Proposition \ref{?impor}, we have an injection $G_{\mathcal{B}, \alpha, \e} \hookrightarrow H_{\mathcal{B}, \alpha, \e}$. Thus, \eqref{p} gives
\begin{equation} \label{biggest}
    N^*_{\alpha,\e}(\mathcal{B}) \ll \#H_{\mathcal{B}, \alpha, \e}.
\end{equation}

The fifth condition of \eqref{paramm} implies that $|X_1^3|, |Y_2^2| \ll \mathcal{B}^{1/6+6\alpha}$, and hence that $|X_1| \ll \mathcal{B}^{1/18+2\alpha}$ and $|Y_2| \ll \mathcal{B}^{1/12+3\alpha}$. The first condition of \eqref{paramm} implies that $|C^3B_2^2D_2| \leq \mathcal{B}$, and hence that $|D_2| \leq \mathcal{B}$. Thus, after forgetting all divisor conditions, \eqref{biggest} becomes
\begin{equation} \label{big}
    N^*_{\alpha,\e}(\mathcal{B}) \ll \underset{|B_2X_1^3|, |CY_2^2| \ll \mathcal{B}^{1/6+6\alpha}}{\underset{|C|, |B_2| \leq \mathcal{B}^{\e}}{\underset{|Y_2| \ll \mathcal{B}^{1/12+3\alpha}}{\sum_{|X_1| \ll \mathcal{B}^{1/18+2\alpha}}}}}\#\left\{(B_1, D_2) \in \Z_{>0} \times \Z_{\neq 0} :\begin{aligned} & |D_2| \leq \mathcal{B} \\ & CY_2^2 - B_2X_1^3 - D_2B_1^6 = 0 \end{aligned} \right\}.
\end{equation}

Now, consider $CY_2^2-B_2X_1^3-D_2B_1^6$ as a polynomial in the variables $B_1$ and $D_2$, with $C, B_2, X_1$, and $Y_1$ fixed. Homogenize this polynomial with the variable $V$ to obtain
\begin{equation*}
    F_{C, B_2, X_1, Y_1}(B_1, D_2, V) := V^7(CY_2^2-B_2X_1^3) -D_2B_1^6,
\end{equation*}
which is a degree $7$ form in the variables $B_1$, $D_2$, and $V$. Given $(B_1,D_2) = (b_1,d_2)$ counted by the right-hand side of \eqref{big}, consider the tuple $(B_1,D_2,V) = (b_1,d_2,1)$. It is a zero of $F_{C, B_2, X_1, Y_1}$, and thus $|D_2B_1^6| \ll \max\{|B_2X_1^3|, |CY_2^2|\} \ll \mathcal{B}^{1/6+6\alpha}$. Hence, we have that $|D_2| \ll \mathcal{B}^{1/6+6\alpha}$ and $|B_1| \ll \mathcal{B}^{1/36+ \alpha}$. However, since $|D_2| \leq \mathcal{B}$, we obtain that $|D_2| \ll \mathcal{B}^{m(\alpha)}$. Moreover, since $B_1D_2 \neq 0$, we must have that $CY_2^2 - B_2X_1^3 \neq 0$.

With this condition, we can show that $F_{C, B_2, X_1, Y_1}$ satisfies the hypotheses of Theorem \ref{bound}. We have that $F_{C, B_2, X_1, Y_1}$ is singular by observing that $(B_1, D_2, V) = (0, 1, 0)$ is a non-zero singular point of $F$. To prove irreducibility over $\Q$, we use Eisenstein's criterion on $\Z[B_1, D_2][V]$ with the prime ideal $(D_2)$ to see that $F_{C, B_2, X_1, Y_1}$ is irreducible over $\Z$. Since $\Z$ is a UFD, we have that $F$ is irreducible over $\Q$ by a generalized version of Gauss' Lemma (e.g., see \cite[Exercise 3.4]{eisenbud1995commutative}).

In sum, for each set counted by the right-hand side of \eqref{big}, the image of an element in the embedding 
\[\left\{(B_1, D_2) \in \Z_{>0} \times \Z_{\neq 0} :\begin{aligned} & |D_2| \leq \mathcal{B} \\ & CY_2^2 - B_2X_1^3 - D_2B_1^6 = 0 \end{aligned} \right\} 
\hookrightarrow \{(B_1,D_2,V) \in \Z_{>0} \times \Z_{\neq 0} \times \Z\}\]
given by $(b_1,d_2) \mapsto (b_1,d_2,1)$ satisfies the following properties:
\begin{itemize}
    \item $(B_1, D_2, V) = (b_1,d_2,1)$ is a zero of $F_{C, B_2, X_1, Y_1}$.
    \item Given such a zero, we have that $|D_2| \ll \mathcal{B}^{m(\alpha)}$, $|B_1| \ll \mathcal{B}^{1/36+\alpha}$, and $|V| \leq 1$.
    \item $F_{C, B_2, X_1, Y_1}$ satisfies the hypotheses of Theorem \ref{bound}.
\end{itemize}

Let $\mathcal{P}(F_{C, B_2, X_1, Y_1})$ be $1$ if the third statement above is true, and let it be $0$ otherwise. In the notation of Theorem \ref{bound}, let $x_1 = B_1$, $x_2 = D_2$,  $x_3 = V$, and the $\mathcal{B}_i$ be the corresponding bounds in the second statement above. Then, using the above properties, \eqref{big} becomes
\begin{align}
    N^*_{\alpha,\e}(\mathcal{B}) & \notag \ll \underset{|C|, |B_2| \leq \mathcal{B}^{\e}}{\underset{|Y_2| \ll \mathcal{B}^{1/12+3\alpha}}{\sum_{|X_1| \ll \mathcal{B}^{1/18+2\alpha}}}}\#\left\{(B_1, D_2, V) \in \Z_{>0} \times \Z_{\neq 0} \times \Z: \begin{aligned} &  F_{C, B_2, X_1, Y_1}(B_1, D_2, V) = 0  \\ & |D_2| \ll \mathcal{B}^{m(\alpha)}, |B_1| \ll \mathcal{B}^{1/36+\alpha}, |V| \leq 1 \\ & \mathcal{P}(F_{C, B_2, X_1, Y_1}) = 1\end{aligned}\right\} \\ & \label{sus} \ll \underset{|C|, |B_2| \leq \mathcal{B}^{\e}}{\underset{|Y_2| \ll \mathcal{B}^{1/12+3\alpha}}{\sum_{|X_1| \ll \mathcal{B}^{1/18+2\alpha}}}} N(F_{C, B_2, X_1, Y_1}; \mathcal{B}_1, \mathcal{B}_2, \mathcal{B}_3)\mathcal{P}(F_{C, B_2, X_1, Y_1}).
\end{align}
Continuing in the notation of Theorem \ref{bound}, we also have
\begin{equation*}
T \ll \mathcal{B}^{m(\alpha) + 6\alpha + 1/6}, \, \mathcal{B}_1\mathcal{B}_2\mathcal{B}_3 \ll \mathcal{B}^{m(\alpha) + \alpha + 1/36}, \,\text{and } k = 7.
\end{equation*}
Finally, we may apply the estimate of Theorem \ref{bound} to \eqref{sus} to obtain
\begin{align} \label{constant}
    N^*_{\alpha,\e}(\mathcal{B}) \ll_{\e} \underset{|C|, |B_2| \leq \mathcal{B}^{\e}}{\underset{|Y_2| \ll \mathcal{B}^{1/12+3\alpha}}{\sum_{|X_1| \ll \mathcal{B}^{1/18+2\alpha}}}} T^{-k^{-2}}(\mathcal{B}_1\mathcal{B}_2\mathcal{B}_3)^{k^{-1}+\e}\ll  \underset{|C|, |B_2| \leq \mathcal{B}^{\e}}{\underset{|Y_2| \ll \mathcal{B}^{1/12+3\alpha}}{\sum_{|X_1| \ll \mathcal{B}^{1/18+2\alpha}}}} \mathcal{B}^{-\frac{1}{49}(m(\alpha) + 6\alpha + \frac{1}{6})} \mathcal{B}^{\frac{1}{7}(m(\alpha) + \alpha + \frac{1}{36})+e(\alpha)}.
\end{align}
Combining terms in \eqref{constant}, we obtain that
\begin{align*}
    N^*_{\alpha,\e}(\mathcal{B})  \ll_{\e} \underset{|C|, |B_2| \leq \mathcal{B}^{\e}}{\underset{|Y_2| \ll \mathcal{B}^{1/12+3\alpha}}{\sum_{|X_1| \ll \mathcal{B}^{1/18+2\alpha}}}}\mathcal{B}^{\frac{1}{49}(6m(\alpha) + \alpha) + c_1 +e(\alpha)} \ll \mathcal{B}^{\frac{6}{49}(m(\alpha) + 41\alpha) + c_1+c_2 +2\e + e(\alpha)},
\end{align*}
as wished.
\end{proof}
The bound obtained in this proposition will be used to obtain bounds for $N_{\alpha}(\mathcal{B})$ in combination with
\begin{prop} \label{step1}
One has that $N_{\alpha}(\mathcal{B}) = N_{\alpha,\e}(\mathcal{B}) +O(\mathcal{B}^{1-\e})$.
\end{prop}
\begin{proof}
We begin by showing that
\begin{equation} \label{boundd}
    \#\{d \in S_6(\mathcal{B}): \square(d) > \mathcal{B}^{\e}\} \ll \mathcal{B}^{1-\e}.
\end{equation}
We can do this directly by proving the chain of inequalities
\begin{equation*} 
    \#\{d \in S_6(\mathcal{B}): \square(d) > \mathcal{B}^{\e}\} \ll \#\left\{ (r,s) \in \Z_{>0}^2 :\begin{aligned}
       & rs^2 \leq \mathcal{B} \\ & s > \mathcal{B}^{\e}
     \end{aligned}\right \} \ll \mathcal{B}^{1-\e}.
\end{equation*}
The first bound follows from the observation that we may uniquely write $d = \mathrm{sgn}(d)d'\square(d)^2$ such that $\mathrm{sgn}(d) \in \{\pm 1\}$ and $d' \in \Z_{>0}$ square-free, and thus the map
\[d \mapsto (d', \square(d)): \{d \in S_6(\mathcal{B}): \square(d) > \mathcal{B}^{\e}\} \to \left\{ (r,s) \in \Z_{>0}^2 :\begin{aligned}
       & rs^2 \leq \mathcal{B} \\ & s > \mathcal{B}^{\e}
     \end{aligned}\right \}\] 
is $2$-to-$1$.

The second bound follows from the observation that
\begin{equation*}
    \left\{ (r,s) \in \Z_{>0}^2 :\begin{aligned}
       & rs^2 \leq \mathcal{B} \\ & s > \mathcal{B}^{\e}
     \end{aligned}\right \} \ll \sum_{j = \lceil \mathcal{B}^{\e} \rceil}^{\lceil \mathcal{B}^{1/2}\rceil}\frac{\mathcal{B}}{j^2} \ll \mathcal{B}\int_{\mathcal{B}^{\e}}^{\mathcal{B}^{1/2}} \frac{1}{j^2} dj \ll \mathcal{B}^{1-\e}.
\end{equation*}
This completes the proof of \eqref{boundd}.

Now, the integers $d$ that $N_{\alpha}(\mathcal{B})$ counts satisfy either $\square(d) \leq \mathcal{B}^{\e}$ or $\square(d) > \mathcal{B}^{\e}$. Those satisfying $\square(d) \leq \mathcal{B}^{\e}$ are counted by $N_{\alpha, \e}(\mathcal{B})$. Those satisfying $\square(d) > \mathcal{B}^{\e}$ are contained in $\{d \in S_6(\mathcal{B}): \square(d) > \mathcal{B}^{\e}\}$, and thus, by \eqref{boundd}, there are $O(\mathcal{B}^{1-\e})$ such $d$. Adding together both counts yields the proposition.
\end{proof}

A consequence of the previous two propositions is

\begin{cor} \label{200}
Let $\e > 0$. The set of $d \in S_6(\mathcal{B})$ such that 
\begin{equation*}
    \zeta_d > |d|^{\resultall-\e}
\end{equation*}
has natural density $1$ in $S_6(\mathcal{B})$ as $\mathcal{B} \to \infty$.
\end{cor}
\begin{proof}
We begin by noting a chain of implications. First, $N_{\alpha, \e}(\mathcal{B}) \leq N^*_{\alpha,\e}(\mathcal{B})$ by the discussion following \eqref{disc}, and so we have that $N^*_{\alpha, \e}(\mathcal{B}) = o(\mathcal{B})$ implies $N_{\alpha,\e}(\mathcal{B}) = o(\mathcal{B})$. Then, from Proposition \ref{step1}, we have that $N_{\alpha,\e}(\mathcal{B}) = o(\mathcal{B})$ implies $N_{\alpha}(\mathcal{B}) = o(\mathcal{B})$. In sum, we have that $N^*_{\alpha, \e}(\mathcal{B}) = o(\mathcal{B})$ implies $N_{\alpha}(\mathcal{B}) = o(\mathcal{B})$, and so by the discussion following \eqref{lookbackto}, it will suffice to show that $\alpha < \frac{43}{246} - \frac{1}{36} = \frac{217}{1476}$ implies $N^*_{\alpha, \e}(\mathcal{B}) = o(\mathcal{B})$. Now, in the notation of Proposition \ref{step2}, suppose that $\alpha$ satisfies
\begin{equation} \label{f1}
\frac{6}{49}(41\alpha + m(\alpha))+ c_1+c_2 < 1.
\end{equation}
Upon taking $\e$ sufficiently small, we then have that \eqref{f1} implies $N^*_{\alpha, \e}(\mathcal{B}) = o(\mathcal{B})$ by Proposition \ref{step2}. Thus, the proof will be complete if we show that $\alpha < \frac{217}{1476}$ implies \eqref{f1}.

We now have two cases:
\begin{itemize}
    \item Case 1: $\alpha \geq \frac{5}{36}$. Then, we have that $m(\alpha) = 1$, and so \eqref{f1} holds if and only if
\begin{equation*}
\frac{246\alpha}{49} + \frac{6}{49} + c_1 + c_2 < 1.
\end{equation*}
This is equivalent to 
\[\alpha < \frac{49(1-c_1-c_2)-6}{246} = \frac{217}{1476}.\]
Thus, \eqref{f1} holds if and only if $\frac{5}{36} \leq \alpha < \frac{217}{1476}$.

    \item Case 2: $\alpha < \frac{5}{36}$. Then, we have that $m(\alpha) = \frac{1}{6} + 6\alpha$, and so \eqref{f1} holds if and only if
\begin{equation*}
    \frac{282\alpha}{49} + \frac{1}{49} + c_1 + c_2  < 1.
\end{equation*}
This is equivalent to 
\[\alpha < \frac{49(1-c_1-c_2)-1}{282} = \frac{247}{1692}.\]
As this bound is weaker than our assumption that $\alpha < \frac{5}{36}$, we conclude that \eqref{f1} holds if and only if $\alpha < \frac{5}{36}$.
\end{itemize}
Combining the bounds in each case, we have that $\alpha < \frac{217}{1476}$ implies \eqref{f1}. This is what we needed to show, and so the corollary follows.
\end{proof}

\bibliographystyle{elsarticle-num} 
\bibliography{final.bib}
\end{document}